\documentclass[12pt]{article}
\usepackage{geometry}
\usepackage{amssymb}
\usepackage{hyperref}
\usepackage{color}

\usepackage{amsmath,amsthm}
\allowdisplaybreaks

\newif\ifleft
\lefttrue 

\newif\ifpre
\prefalse

\newif\ifhide
\hidefalse

\def\1{\mathbf{1}}

\def\P{\mathsf{P}}

\newtheorem{theorem} {Theorem}
\newtheorem{Proposition}{Proposition}

\newtheorem{remark} {Remark}
\newtheorem{corollary} {Corollary}

\newif \ifshowup
\showupfalse 

\title{{\normalsize\tt\hfill\jobname.tex}\\
      \bf On weak existence of
solutions of degenerate McKean--Vlasov
equations} 
\author{A.Yu. Veretennikov\footnote{Institute for Information Transmission Problems,  Moscow, Russian Federation; email: ayv@iitp.ru.
}
}

\begin{document}

\maketitle


\begin{abstract}
A new weak 
existence result for degenerate multi-dimensional stochastic McKean--Vlasov equation is established  under relaxed regularity conditions. 

~


~


\end{abstract}

\maketitle

\section{Introduction}
The subject of this paper is solutions  of the stochastic It\^o-McKean-Vlasov (McKean-Vlasov) equation in $\mathbb R^{2d}$
\begin{equation}\label{e1}
dX_t = B_0[t,Z_t, \mu_t] dt, \quad 
dY_t = B_1[t,Z_t, \mu_t]dt + \Sigma[t,Z_t, \mu_t]dW_t, \quad
X_0=x_0,\,Y_0=y_0,
\end{equation}
where $Z_t = (X_t, Y_t) \in \mathbb R^{2d}$, in a particular situation called the true McKean-Vlasov case  under the convention
\begin{equation}\label{e200}
B_i[t,z,\mu]=\int
b_i(t,z,\zeta)\mu(d\zeta)\; \; (i=0,1),\;\; \Sigma[t,z,\mu]=\int
\sigma(t,z,\zeta)\mu(d\zeta),
\end{equation}
where $z = (x,y)\in \mathbb R^{2d}$ and $\zeta = (\xi,\eta)\in \mathbb R^{2d}$,  
and under certain non-degeneracy assumptions on $\sigma$.
Here $W$ is a standard $d$-dimensional
Wiener process, $b$ and $\sigma$ are vector and matrix Borel
functions of corresponding dimensions $d$ and $d\times
d$, $\mu_t$ is the distribution of the process $Z$ at time
$t$. The initial data $z_0=(x_0,y_0)$ may be random and in this case it is  independent of $W$. Vlasov's proposal was a substitution of a real multiparticle  interaction by a certain ``mean field'' \cite{Vlasov68}. The classical introduction to the whole topic in its stochastic version may be found in \cite{Sz}; one more very important  reference  is \cite{Dobr}, although it is  devoted soleyly  to the deterministic setting. In the present paper we investigate only the problem of weak existence for genuinly degenerate stochastic McKean -- Vlasov equations. The equations like (\ref{e1}) naturally arise in mechanical systems with stochastic forces or noise. The aim of this paper is to show weak existence for a such a degenerate SDE system under minimal regularity assumptions on both coefficients with respect to all variables. We assume the non-degeneracy of $\sigma$ and highlight that this non-degeneracy only holds for the second component $Y$ of the system (\ref{e1}). The interest to the minimal regularity is mainly due to the control problems where the optimal strategies are usually discontinuous.  

\medskip

{ 
Among important works on the subject including more recent ones, there are the papers \cite{Funaki, BossyJabir,  CarmonaDelarue1, CarmonaDelarue2, Issoglio, Mehri,  RoecknerZhang, TalayTom, Tom, Xicheng}; see also the bibliography therein. Let us mention a new extended setting in the publications  
\cite{Xicheng, TalayTom} where coefficients of the equation may depend on the marginal density of the process beside a more ``usual'' dependence on the marginal distribution. In this paper we do not pursue this goal. If in the first equation in (\ref{e1}) the function $B_0$ were equal identically $Y_t$, then our results would be close to those in a special case ``without $\rho$'' in  \cite{Xicheng}, 
except that the form of the coefficient $\Sigma$ in \cite{Xicheng} is different from ours and that we do not consider unbounded drifts. In particular, our $B_0$ is also bounded, which, of course, may be relaxed. 
Our general setting as in (\ref{e1}) apparently is not covered by results in \cite{Xicheng}. Also, we do not touch here the issue of weak or strong uniqueness. 
}

\medskip


\medskip


The study in the present paper is based on Krylov's bounds \cite{Kry} on Skorokhod's technique of weak convergence, and on the approach proposed in \cite{Nisio} for the ordinary It\^o SDEs which was further generalised to some extent in \cite{Ver-ait} also for the ordinary It\^o SDEs. Other useful references may be found in the cited papers.

The structure of the paper is as follows.   In the section \ref{sec:we}  weak existence is stated under appropriate conditions. The section 3 containts its proof based on a combination of Krylov's bounds and on Krylov's existence results for the nondegenerate Ito's equations \cite{Kry69}, \cite{Kry}, and on Nisio's  weak existence proof also for Ito's SDEs \cite{Nisio}. The degeneracy of the diffusion is overcome still by using Krylov's bounds for {\em non-degenerate} It\^o processes. No regularity  of the coefficients $b_1$ and $\sigma$ is assumed with respect to the  variables $t$, $y$, and $\eta$. Uniform continuity is assumed for both $b_1$ and $\sigma$ with respect to the  variables $x$ and $\xi$, and for $b_0$ with respect to  all variables except $t$. 
The abbreviation by CBS signifies the Cauchy-Buniakovsky-Schwarz inequality and BCM stands the Bienaym\'e-Chebyshev-Markov inequality.

\section{Weak existence}\label{sec:we}
\subsection{Main results}
Let us recall a fact
from functional analysis useful for the case  (\ref{e1})-- (\ref{e200}),
see, for example, (see \cite[Theorem 1.5.5]{Kry-ln}). The proposition \ref{Pro1} and its corollary are stated in a slightly more general form than what is needed for bounded coefficients. 

\begin{Proposition}\label{Pro1}
For any Borel function $f(z, \zeta)$ and any probability measure
$\mu(d\zeta)$ such that $f(z, \cdot)$ is integrable with respect to this measure, the
function $\displaystyle f[z, \mu]:= \int f(z,\zeta)\,\mu(d\zeta)$ is   Borel measurable in $z$.
\end{Proposition}

\begin{corollary}\label{cor1}
Suppose for each $(t,z)$ the Borel coefficients $b(t,z,\zeta)$ and $\sigma(t,z,\zeta)$ are bounded in $\zeta$ and integrable in \(z\) with respect to all $(\mu_t), \, t\ge 0$, where \(\mu_t\) are marginal distributions of any weak solution of the equation (\ref{e1}).  Then the functions $\tilde b(t,z):=B[t, z, \mu_t]$ and $\tilde \sigma(t,z):=
\Sigma[t, z, \mu_t]$ are Borel measurable in $(t, z)$.
\end{corollary}

Recall the notations $z=(x,y)$, $\zeta = (\xi, \eta)$.
\begin{theorem}\label{thm1}
Let the initial value \(z_0\) have a finite fourth moment and assume that  the following three conditions are satisfied. 
(1$^\circ$) Firstly, the functions $b_i, \, i=0,1$, and $\sigma$ are uniformly bounded, i.e.,
there exists $C>0$ such that for any { $s,z,\zeta$},
\begin{equation}\label{linear}
|b_0(s,z,\zeta)|+|b_1(s,z,\zeta)|+\|\sigma(s,z,\zeta)\|\leq C,
\end{equation}
where $|\cdot|$ stands for the Euclidean norm in $\mathbb R^d$ for
$b_i$ and $\|\cdot \|$ for the  $\|\sigma\| =
\sqrt{\sum_{i,j}\sigma_{ij}^2}\,$. 
(2$^\circ$) Secondly, the diffusion matrix $\sigma(s,z,\zeta)$ is symmetric and  
uniformly nondegenerate in the following sense: there is a value $\nu>0$ such that 
\begin{equation}\label{si}
\inf\limits_{s,z,\zeta}\inf\limits_{|\lambda|=1}
\lambda^*\sigma(s,z,\zeta)
\lambda 
\ge \nu.
\end{equation}
(3$^\circ$) Thirdly, $b_1(t,x,y,\xi,\eta)$ and $\sigma(t,x,y,\xi,\eta)$ are 
continuous with respect to $(x,\xi)$ for each $(t,y,\eta)$ with a uniform modulus of continuity $\rho(\cdot)$. 
\\ 
(4$^\circ$) Finally, for any $t$ the function $b_0(t, z, \zeta)$ is continuous in the variables $(z, \zeta)$.

Then the equation (\ref{e1}) has a weak solution on some probability space with a standard $d$-dimensional
Wiener process with respect to some filtration $({\cal F}_t, \, t\ge 0)$. 

\end{theorem}

\begin{remark}
It is likely that 
the requirement of the fourth moment of $z_0$ may be relaxed. Also, it is likely that the condition of boundedness of the drift itself may be relaxed considerably. We do not pursue both goals in order to simplify the presentation and references. 
\end{remark}

\medskip

Denote 
$$
A[t,z,\mu]:= \Sigma\Sigma^*[t,z,\mu].
$$

\subsection{Proof}
{\bf 1.} 
Let us  mollify all three coefficients $b_0, b_1$,  and $\sigma$ with respect to all variables by convolutions in such a way that they become globally Lipschitz in $z$,  $\zeta$, and $t$. Namely, let for $i=0,1$,
\begin{equation*}
b_i^n(t,z,\zeta) = b(t,z,\zeta) * \psi_n(t) * \varphi_n(x) * \varphi_n(y) * 
\varphi_n(\xi) * \varphi_n(\eta),
\end{equation*}
and
\begin{equation*}
\sigma^n(t,z,\zeta) = \sigma(t,z,\zeta) * \psi_n(t) * \varphi_n(x) * \varphi_n(y) * \varphi_n(\xi) * \varphi_n(\eta),
\end{equation*}
 where the sequences $\varphi_n(\cdot)$ and $\psi_n$ are defined in a
standard way, i.e., as non-negative \(C^\infty\) functions with a compact support, integrated to one, and so that this compact support squeezes to the origin of the corresponding variable as \(n\to\infty\); or, in other words, that they are delta-sequences in the corresponding variables. {For example, denoting $z'=(x',y')$ and $\zeta' = (\xi', \eta')$,
\begin{align*}
&b_i(t,z,\zeta) * \psi_n(t) * \varphi_n(x) * \varphi_n(y) * 
\varphi_n(\xi) * \varphi_n(\eta) 
 \\\\
&:= \!\int b_i(t',z',\zeta')\psi_n(t\!-\!t') \varphi_n(x\!-\!x') \varphi_n(y\!-\!y') 
\varphi_n(\xi\!-\!\xi') \varphi_n(\eta\!-\!\eta')\,dt' dx' dy' d\xi' d\eta'.
\end{align*}
The (multiple) integral without limits always means integration over the whole domain, in our case over $\mathbb R^{1+4d}$. }
Note that, of course, for every $n$ the smoothed coefficients remain uniformly bounded and all have the same uniform moduli of continuity with respect to the variables $y,\eta$; also, the smoothed diffusion $\sigma$ remains uniformly non-degenerate with ellipticity constants independent of~$n$. 
While performing the convolution with $\psi_n$, it is assumed that $\sigma(t,z,\zeta)\equiv I_{d\times d}$ for $t<0$ (this is needed to leave the mollified diffusion acting on the variable $y$ uniformly nondegenerate for $t\ge 0$ near zero), and that $b_i(t,z,\zeta)\equiv 0$ for $t<0$ and $i=0,1$.

\medskip

The equation with smoothed coefficients has a strong solution. Even under weaker linear growth conditions it is explained, for example, in \cite[proof of proposition 1]{VerMish20}, as well as in many other sources; this is not linked to the non-degeneracy in any way.

\medskip

\noindent
{\bf 2}. In a standard way (see, e.g.,   the proof of \cite[theorem 1.6.4]{GS68}), the estimates uniform in $n$ follow,
\begin{align}
\label{kol1}
 \mathsf{E} \sup_{0\le t\le T} |Z_t^n |^4 \le C_T (1+\mathsf E |z_0|^4),
\end{align}
and
\begin{eqnarray}\label{kol2}
\sup_{0\le s\le t\le T; \, t-s\le h}\mathsf{E} |Z_t^n - Z^n_s|^4 \le
C_{T}  h^2,
\end{eqnarray}
with some constants $C_{T}$ which may be different for different inequalities but do not depend on $n$. (In fact, in \cite{GS68} the assumptions allow a linear growth in $x$; overall, it is a very standard material.)

\medskip

\noindent
{\bf 3}.
 Let us introduce new processes $(\xi^n,\eta^n)=:\zeta^n$ which are the copies of  $(X^n,Y^n)=:Z^n$, that  satisfy similar SDEs on some independent  probability spaces. In the sequel by $ \mathsf{E}^3  \sigma^n(s,{ Z^n_s,\zeta^n_s}) $ 
we denote expectation with respect to the third variable $ \zeta^n_s $ 
{\em conditional} on { $Z^n_s$}, 
that is, 
$$ 
\mathsf{E}^3  \sigma^n(s,Z^n_s,\zeta^n_s) = \int \sigma^n(s,Z^n_s,\zeta)\mu^{\zeta^n}_s(d\zeta),
$$  where $ \mu^{\zeta^n}_s = {\cal L}(\zeta^n_s)$ (here naturally $\zeta$ is the variable of integration). {  Equivalently, it may be written as $\mathsf E \sigma(s,z,\zeta^n_s)\vert_{z=Z^n_s}$.} Likewise,
$$
\mathsf{E}^3  (\sigma^n(s,Z^n_s,\zeta^n_s)
- \sigma^n(s,Z^0_s,\zeta^0_s))
$$
is another notation for 
$$
\displaystyle \int \sigma^n(s,Z^n_s,\zeta)\mu^{\zeta^n}_s(d\zeta)
-  \int \sigma^n(s,Z_s,\zeta)\mu^{\zeta^0}_s(d\zeta),
$$
where  $\mu^{\zeta}_s = {\cal L}(\zeta_s)$ for any random variable $\zeta \in \mathbb R^{2d}$; 
the integral 
$$
\mathsf{E}^3 \|\sigma^n(s, Z^n_s, { \zeta^n_s}) - \sigma(s,Z^0_s, { \zeta^0_s})\|^2
$$
is understood as
$$
\displaystyle \int \|\sigma^n(s,Z^n_s,\zeta)
-  \sigma^n(s,Z_s,\zeta')\|^2
{\mathsf P}({ \zeta^n_s} \in d\zeta, { \zeta^0_s}\in d\zeta'),
$$
if $\zeta^n$ and $\zeta^0$ are defined on the same probability space.

\medskip

Due to the estimates  (\ref{kol1})--(\ref{kol2}) and by virtue of Skorokhod's
Lemma about a single probability
space and convergence in probability (see \cite[\S 6, ch. 1]{Sko}, or \cite[Lemma 2.6.2]{Kry}, or \cite[Lemma 4 in the Appendix]{VerMish20}) 
without loss of generality we may and will  assume that
not only $\mu^n \Longrightarrow \mu$, but also on some probability space
\begin{equation}\label{converge}
(\tilde Z^{n'}_{t},\tilde \zeta^{n'}_{t}, \tilde W^{n'}_{t}) \stackrel{\mathsf{P}}{\to}
(\tilde Z^0_t,\tilde \zeta^0_t, \tilde W^0_t), \quad n\to\infty,
\end{equation}
{ generally speaking, over a sub-sequence $n'\to\infty$, for any~$t$ and  for some   equivalent  random processes $(\tilde Z^{n'}, \tilde \zeta^{n'}, \tilde W^{n'})$, and redenote the subsequence $(n')$ again by $(n)$.}

Slightly abusing notations, we  denote initial values still by \(z_0\) without tilde.
Also,  without loss of generality we assume that each process $(\tilde \zeta^n_{t}, \, t\ge 0)$ for any $n\ge 1$ is independent of $(\tilde Z^n, \tilde W^n)$, as well as their limit $\tilde \zeta^0_t$ may be chosen to be independent of the limits $(\tilde Z^0,\tilde W^0)$ (this follows from the fact that on the original probability space $\eta^n$ is independent of $(Z^n, W^n)$ and on the new probability space their joint distribution remains the same; hence, independence of $\tilde \zeta^n$ is also valid and in the limit this is still true). See the details in the proof of the Theorem 2.6.1 in \cite{Kry}. 
On  independent probability spaces we have,
\begin{equation}\label{exi0}
d\xi^n_t = B_0[t,\zeta^n_t,\mu_t]dt, \; d\eta^n_t = B_1^{n}[t,\zeta^n_t, \mu_t]dt + \Sigma^{n} [t,\zeta^n_t, \mu_t]dW^{\prime,n}_t, \; t\ge 0, \;
{\cal L}(\zeta^n_0)={\cal L}(z_0),
\end{equation}
and
\begin{equation*}
d\tilde \xi^n_t = B_0[t,\tilde \zeta^n_t,\mu_t]dt, \; d\tilde \eta^n_t = B^{n}[t,\tilde \zeta^n_t, \mu_t]dt +  \Sigma^{n} [t,\tilde \zeta^n_t, \mu_t]d\tilde W^{\prime,n}_t, \;\; t\ge 0, \;
{\cal L}(\tilde \zeta^n_0)={\cal L}(z_0).
\end{equation*}

Due to the inequality (\ref{kol2}),
the same inequality holds for $\tilde Z^n$ and  $\tilde W^n$, in particular,
\begin{eqnarray}\label{kol2t}
\sup_{0\le s\le t\le T; \, t-s\le h}\mathsf{E} |\tilde Z_t^n - \tilde Z^n_s|^4 \le C_{T}  h^2.
\end{eqnarray}
Due to Kolmogorov's continuity theorem, it means that all processes $\tilde Z^n$ may be regarded as continuous, and  $\tilde W^n$ can be also assumed continuous by the same reason. Note for the sequel that the bound (\ref{kol1}) is also applicable to the process $\tilde Z$:
\begin{align}
\label{kol11}
 \mathsf{E} \sup_{0\le t\le T} |\tilde Z_t^n |^4 \le C_T (1+\mathsf E |z_0|^4)
\end{align}
because of the equivalence of $Z$ and $\tilde Z$. 

Further, due to the independence of the increments of $W^n$ after time $t$ of the sigma-algebra $\sigma(Z^n_s, W^n_s, s\le t)$, the same property holds true for $\tilde W^n$ and  $\sigma(\tilde Z^n_s, \tilde W^n_s, s\le t)$, as well as for $\tilde W^n$ and for the completions of the sigma-algebras  $\sigma(\tilde Z^n_s, \tilde W^n_s, s\le t)$ which we denote by ${\cal F}^{(n)}_t$. Also,  the processes  $\tilde Z^n$ are adapted to the filtration $({\cal F}^{(n)}_t)$. So, all stochastic integrals which involve $\tilde Z^n$ and $\tilde  W^n$ are well defined. The same relates to the processes $\tilde \zeta^n$.

So,  we may now hope to pass to the limit { as $n\to\infty$} in the equation
\begin{align*}
&\tilde X^{n}_t = x_0 + \int_0^t \mathsf{E}^3 b_0^{n}(s,\tilde Z^{n}_s, \tilde \zeta^{n}_s)\,ds, 
 \\\\
&\tilde Y^{n}_t = y_0 + \int_0^t \mathsf{E}^3 b_1^{n}(s,\tilde Z^{n}_s, \tilde \zeta^{n}_s)\,ds +
\int_0^t \mathsf{E}^3 \sigma^{n}(s,\tilde Z^{n}_s, \tilde \zeta^{n}_s) d\tilde W^{n}_s,
\end{align*}
in order to get
\begin{align*}
&\tilde X^{0}_t = x_0 + \int_0^t \mathsf{E}^3 b_0^{}(s,\tilde Z^{0}_s, \tilde \zeta^{0}_s)\,ds, 
 \nonumber \\\\  \nonumber 
&\tilde Y^{0}_t = y_0 + \int_0^t \mathsf{E}^3 b_1^{}(s,\tilde Z^{0}_s, \tilde \zeta^{0}_s)\,ds +
\int_0^t \mathsf{E}^3 \sigma^{}(s,\tilde Z^{0}_s, \tilde \zeta^{0}_s) d\tilde W^{0}_s,
\end{align*}
or, equivalently,
\begin{align}\label{limitsde}
&\tilde X^{0}_t = x_0 + \int_0^t B_0^{}(s,\tilde Z^{0}_s, \mu_s)\,ds, 
 \nonumber \\\\  \nonumber 
&\tilde Y^{0}_t = y_0 + \int_0^t B_1^{}(s,\tilde Z^{0}_s, \mu_s)\,ds +
\int_0^t \Sigma^{}(s,\tilde Z^{0}_s, \mu_s) d\tilde W^{0}_s, \quad \mu_s = {\cal L}(\tilde Z^0_s).
\end{align}

Recall that a priori bounds
(\ref{kol1}) -- (\ref{kol2}) and (\ref{kol2t}) hold true with constants not depending on $n$. 
By virtue of the  a priori estimates for \(\tilde W^n\), the process \(\tilde W^0\) is continuous and it is, naturally, a $d$-dimensional Wiener process. By virtue of the uniform estimates (\ref{kol2}), the limit \((\tilde Z^{0}_t,\tilde  \zeta^{0}_t)\) may also be regarded as  continuous due to Kolmogorov's continuity theorem,  because   a priori bounds  (\ref{kol1}) -- (\ref{kol2}) remain valid for the limiting processes \(\tilde Z, \tilde \zeta\). 

\medskip

\noindent
{\bf 4.}
Before going further we must make sure that the stochastic integral $\int_0^t \mathsf{E}^3 \sigma^{}(s,\tilde Z^{0}_s, \tilde \zeta^{0}_s) d\tilde W^{0}_s$ is well-defined. Denote ${\tilde {\cal F}^0_t} = \sigma(\tilde Z^0_s, \tilde W^0_s, \, s\le t)$. For our goal it suffices to prove that $(\tilde W^0_t,{\tilde {\cal F}^0_t})$ is a Wiener process, that is, that the increments $\tilde W^0_t - \tilde W^0_s$ are independent of the sigma-algebra ${\tilde {\cal F}^0_s}$ for any $t>s$. In turn, for this aim it suffices to show that  for any $k$ and any $s_1 < s_2 <\ldots <s_k<s_{k+1}$ and any compacts $A_i \in {\cal B}(\mathbb R^{2d}), \, i=1, \ldots, k$ and $B_i \in {\cal B}(\mathbb R^{d}), \, i=1, \ldots, k+1$ the following two probabilities are equal: 
\begin{align}\label{inde}
& \mathsf P\left(\bigcap_{i=1}^k (\tilde Z^0_{s_i}\in A_i) \bigcap \bigcap_{i=1}^k (\tilde W^0_{s_i}\in B_i)) \bigcap (\tilde W^0_{s_{k+1}} - \tilde W^0_{s_{k}}\in B_{k+1})\right) 
\nonumber
  \\\\
\nonumber
& = \mathsf P\left(\bigcap_{i=1}^k (\tilde Z^0_{s_i}\in A_i) \bigcap \bigcap_{i=1}^k (\tilde W^0_{s_i}\in B_i)\right) \mathsf P \left(\tilde W^0_{s_{k+1}} - \tilde W^0_{s_{k}}\in B_{k+1}\right).
\end{align} 
Consider the cylinders
$$
C_1:= \left(\bigcap_{i=1}^k (z_i \in A_i) \bigcap \bigcap_{i=1}^k ( w_i\in B_i) \bigcap (v\in B_{k+1})\right) \subset \mathbb R^{(2k+1)d},  
$$
$$
C_2:= \left(\bigcap_{i=1}^k (z_i \in A_i) \bigcap \bigcap_{i=1}^k (w_i\in B_i)\right) \subset \mathbb R^{2kd},  
$$
and
$$
C_3:= (v\in B_{k+1}) \subset \mathbb R^{d}.  
$$
Let us introduce the following measures:
$$
\nu_1(C_1) := \mathsf P \left((\tilde Z^0_{s_1}, \ldots, \tilde Z^0_{s_k}, \tilde W^0_{s_1}, \ldots, \tilde W^0_{s_k}, \tilde W^0_{s_{k+1}} - \tilde W^0_{s_{k}}) \in C_1\right), 
$$
$$
\nu_2(C_2) := \mathsf P \left((\tilde Z^0_{s_1}, \ldots, \tilde Z^0_{s_k}, \tilde W^0_{s_1}, \ldots, \tilde W^0_{s_k}) \in C_2\right), 
$$
and
$$
\nu_3(C_3) := \mathsf P \left(\tilde W^0_{s_{k+1}} - \tilde W^0_{s_{k}} \in C_3\right). 
$$
All three measures are naturally uniquely extended to Borel sigma-algebras, respectively, in $\mathbb R^{2k+1}$,  $\mathbb R^{2k}$, and  $\mathbb R$.
Now, we want to justify the equality 
\begin{equation}\label{nu123}
\nu_1(C_1) = \nu_2(C_2)\nu_3(C_3),
\end{equation}
which is equivalent to the desired equation (\ref{inde}). As it is known (see, for example, \cite[Theorem 1.2.4]{Kry-ln}), any sigma-finite measure on $\mathbb R^d$ is uniquely determined by its integrals with all bounded continuous functions. Hence, it suffices to establish the following equality for any two functions $f\in C_b(\mathbb R^{2k})$ and $g\in C_b(\mathbb R)$, 
\begin{align}\label{nu123}
&\int f(z_1, \ldots, z_k, w_1, \ldots, w_k) g(w_{k+1})\nu_1(dz_1 \ldots dz_k dw_1 \ldots dw_k dv) 
\nonumber
 \\\\
 \nonumber
&= \int f(z_1, \ldots, z_k, w_1, \ldots, w_k) \nu_2(dz_1 \ldots dz_k dw_1 \ldots dw_k) \int g(v)\nu_3(dv) .
\end{align}
The latter equality (\ref{nu123}) may be rewritten in the form
\begin{align}\label{nu1230}
&\mathsf E f(\tilde Z^0_{s_1}, \ldots, \tilde Z^0_{s_k}, \tilde W^0_{s_1}, \ldots, \tilde W^0_{s_k}) g(\tilde W^0_{s_{k+1}} - \tilde W^0_{s_{n}})
\nonumber
 \\\\
 \nonumber
&= \mathsf E f(\tilde Z^0_{s_1}, \ldots, \tilde Z^0_{s_k}, \tilde W^0_{s_1}, \ldots, \tilde W^0_{s_k}) \mathsf E g(\tilde W^0_{s_{k+1}} - \tilde W^0_{s_{k}}).
\end{align}
But this equality for continuous bounded $f$ and $g$ immediately follows from the pre-limiting equation which is valid for each $n>0$:
\begin{align}\label{nu1230}
&\mathsf E f(\tilde Z^n_{s_1}, \ldots, \tilde Z^n_{s_k}, \tilde W^n_{s_1}, \ldots, \tilde W^n_{s_k}) g(\tilde W^n_{s_{k+1}} - \tilde W^n_{s_{k}})
\nonumber
 \\\\
 \nonumber
&= \mathsf E f(\tilde Z^n_{s_1}, \ldots, \tilde Z^n_{s_k}, \tilde W^n_{s_1}, \ldots, \tilde W^n_{s_k}) \mathsf E g(\tilde W^n_{s_{k+1}} - \tilde W^n_{s_{k}}).
\end{align}
Indeed, we know that solution $\tilde Z^n$ is strong, so the latter equality is the corollary of the fact that $\tilde W^n$ is a Wiener process with respect to its own family of sigma-fields $({\cal F}^{\tilde W^n}_t)$ which coincides with $({\cal F}^{\tilde W^n, \tilde Z^n}_t)$. Hence, indeed, the stochastic integral $\int_0^t \mathsf{E}^3 \sigma^{}(s,\tilde Z^{0}_s, \tilde \zeta^{0}_s) d\tilde W^{0}_s$ is well-defined, as  promised.

\medskip

\noindent
{\bf 5}.
\noindent
We have to show that
\begin{equation}\label{limNb0}
\int_0^t \mathsf{E}^3 b_i^{n}(s,\tilde Z^{n}_s, \tilde \zeta^{n}_s)ds  \stackrel{\mathsf P}{\to} \int_0^t \mathsf{E}^3 b_i(s,\tilde Z^{0}_s, \tilde \zeta^{0}_s)ds, \;\; i=0,1, 
\end{equation}
and
\begin{equation}\label{limNs0}
\int_0^t \mathsf{E}^3 \sigma^{n}(s,\tilde Z^{n}_s,  \tilde \zeta^{n}_s) d\tilde W^{n}_s  \stackrel{\mathsf P}{\to} \int_0^t  \mathsf{E}^3 \sigma(s,\tilde Z^{0}_s, \tilde  \zeta^{0}_s) d\tilde W^{0}_s, \quad n\to\infty.
\end{equation}
We start with the drift term. The cases $b_0$ and $b_1$ are tackled quite similarly, so let us show one of them for $b_0$. 

\medskip

{  Let $c>0$.} 
Let us fix some  $n_0$ and let $n>n_0$; { in the sequel both $n$ and $n_0$ will tend to infinity}. 
We have for any $t\le T$,
\begin{align*}
& \displaystyle \mathsf P \left(\left|\int_0^t \left( (\mathsf{E}^3  b_0^{n}(s,\tilde X^{n}_s, \tilde Y^{n}_s, \tilde \xi^{n}_s, \tilde \eta^{n}_s))  -  (\mathsf{E}^3  b_0(s,\tilde X^0_s,\tilde Y^0_s, \tilde \xi^{0}_s,\tilde \eta^{0}_s))\right)ds\right| > c\right)
 \nonumber \\\nonumber \\\nonumber
& \!\!\!\displaystyle \le \!\mathsf P\!\left(
\left|\int_0^t \!\!\! \left(\!\mathsf{E}^3 b_0^{n}(s,\tilde X^{n}_s, \tilde Y^{n}_s, \tilde \xi^{n}_s,\tilde \eta^{n}_s) \! -\! (\mathsf{E}^3 b_0^{n_0}(s,\tilde X^{n}_s, \tilde Y^{n}_s, \tilde \xi^{n}_s, \tilde \eta^{n}_s)\right)ds\right| \!>\! \frac{c}{3}\right)
 \nonumber \\\nonumber \\\nonumber
& \displaystyle + \mathsf P \left(\left|\int_0^t  \left( ( \mathsf{E}^3 b_0^{n_0}(s,\tilde Z^{n}_s, \tilde \zeta^{n}_s))  -  ( \mathsf{E}^3 b_0^{n_0}(s,\tilde Z^0_s,\tilde  \zeta^0_s))\right)ds\right| > \frac{c}3\right)
 \nonumber \\\nonumber \\\nonumber
& \!\displaystyle \!+\! \mathsf P \left(\left|\int_0^t \!\!\! \left(\mathsf{E}^3 b_0^{n_0}(s,\tilde X^{0}_s, \tilde Y^{0}_s, \tilde \xi^{0}_s,\tilde \xi^{0}_s))) \!-\!  \mathsf{E}^3 b_0(s,\tilde X^{0}_s, \tilde Y^{0}_s, \tilde \xi^{0}_s, \tilde \xi^{0}_s))\right)ds\right| \!>\! \frac{c}3\right)
  \\ \nonumber\\
& \displaystyle =: I^1 + I^2 + I^3.
\end{align*}

Now the idea is that on a finite interval of time on each $\omega$ the components $\tilde X^{n}_s$ and $\tilde \xi^{n}_s$ are close to certain trajectories of some countable epsilon-net of continuous (even differentiable) functions in $C([0,T]; \mathsf R^d)$. { (Notice that for the components $\tilde Y^{n}_s$ and $\tilde \eta^{n}_s$ the same is also true, but we will not use it in what follows.)} Denote this net by ${\cal N}_\epsilon$ and the union of its first $N$ elements by ${\cal N}_{N,\epsilon}$. More than that, since $\tilde Z^{n}_s$ and $\tilde \zeta^{n}_s$ are bounded in probability (uniformly in $n$) on any bounded interval $[0,T]$, we may take into account only finitely many elements of this epsilon-net, up to a small enough probability, that is, for any $\epsilon >0$ there exists $M>0$ such that 
{ 
\begin{equation}\label{eps1}
\sup_{n\ge 1} \P(\sup_{0\le t\le T}|\tilde Z^n_t| \vee |\tilde \zeta^n_t| > M)<\epsilon,
\end{equation}
and there exists $N>0$ such that 
\begin{equation}\label{eps2}
\sup_{n\ge 1}\,\P(\bigcup_{k,j=1}^N\sup_{0\le t\le T}|\tilde X^n_t - \phi^k_t| \vee |\tilde \xi^n_t - \phi^j_t| > \epsilon)<\epsilon,
\end{equation}
}
where all $\phi^k, \phi^j \in {\cal N}_\epsilon$. The value $N$ may be chosen uniformly with respect to $n$  due to the a priori bounds (\ref{kol2t})--(\ref{kol11}) and, moreover, because the trajectories $(\tilde X^n_t)$  and $(\tilde \xi^n_t)$  admit { Lipschitz bounds independent on $n$.} 

\medskip

{ 
Let
$$
A^n_{M,\epsilon}:= (\omega: \sup_{0\le t\le T}|\tilde Z^n_t| \vee |\tilde \zeta^n_t| > M),
$$
and 
$$
B^n_{N,\epsilon}:= (\omega: 
\bigcup_{k,j=1}^N\sup_{0\le t\le T}|\tilde X^n_t - \phi^k_t| \vee |\tilde \xi^n_t - \phi^j_t| > \epsilon).
$$
}
Outside these two events $A^n_{M,\epsilon}$ and $B^n_{N,\epsilon}$ of the total probability not exceeding $2\epsilon$ we may assume that 
$$
\sup_{0\le t\le T}|\tilde Z^n_t| \vee |\tilde \zeta^n_t| \le  M, 
$$
and 
\begin{equation}\label{phiphi}
\inf_{k,j\le N} \sup_{0\le t\le T}|\tilde X^n_t - \phi^k_t| \vee |\tilde \xi^n_t - \phi^j_t| \le \epsilon.
\end{equation}
On the event $A^n_{k,j,\epsilon}:=(
\sup_{0\le t\le T}|{ \tilde X^n_t} - \phi^k_t| \vee |\tilde { \xi^n_t} - \phi^j_t| \le \epsilon)$ we have, 
\begin{align*}
\left|\int_0^t  \!\!\mathsf{E}^3  1(A^n_{k,j,\epsilon}) \left(b_0^{n}(s,\tilde X^{n}_s, \tilde Y^{n}_s, \tilde \xi^{n}_s,\tilde \eta^{n}_s) \! -\!  b_0^{n}(s,\phi^k_s, \tilde Y^{n}_s,\phi^j_s, \tilde \eta^{n}_s)\right)ds\right| 
\!\le \! t\rho(\epsilon).
\end{align*}
where $\rho$ is the joint modulus of continuity of both coefficients $b(s,x,y,\xi,\eta)$ and $\sigma(s,x,y,\xi,\eta)$ in $x$ and in $\xi$.

\medskip

Similar bounds hold true for the pair $(\tilde X^0_t, \tilde \xi^0_t)$ due to the convergence and because of the a priori bounds (\ref{kol1}). Therefore, there exists $M$  such that 
$$
\P(\underbrace{\sup_{0\le t\le T}|\tilde Z^0_t| \vee |\tilde \zeta^0_t| > M }_{=:A^0_{M,\epsilon}})<\epsilon,
$$
and there exists $N>0$ such that 
$$
\P(\underbrace{\bigcup_{k,j=1}^N\sup_{0\le t\le T}|\tilde X^0_t - \phi^k_t| \vee |\tilde \xi^0_t - \phi^j_t| > \epsilon }_{=:B^0_{N,\epsilon}})<\epsilon,
$$
where all $\phi^k, \phi^j \in {\cal N}_{N,\epsilon}$.

Replacing $\tilde X^{n}_s$ and $\tilde \xi^{n}_s$ by nonrandom $\phi, \psi\in {\cal N}_{N,\epsilon}$ in the integrals like 
\begin{equation}\label{xphipsin}
\mathsf P\left(
\left|\int_0^t  \left(\mathsf{E}^3 b_0^{n}(s,\phi_s, \tilde Y^{n}_s, \psi_s,\tilde \eta^{n}_s)  -  (\mathsf{E}^3 b_0^{n_0}(s,\phi_s, \tilde Y^{n}_s, \psi_s, \tilde \eta^{n}_s)\right)ds\right| > \frac{c}{3}\right), 
\end{equation}
we will be able to apply Krylov's bounds to show convergence due to the nondegeneracy of $\sigma$; a similar approach is applicable to the probability 
\begin{equation}\label{xphipsi}
\mathsf P\left(
\left|\int_0^t  \left(\mathsf{E}^3 b_0^{n_0}(s,\phi_s, \tilde Y^{0}_s, \psi_s, \tilde \xi^{0}_s)  -  (\mathsf{E}^3 b_0^{}(s,\phi_s, \tilde Y^{0}_s, \psi_s, \tilde \xi^{0}_s)\right)ds\right| > \frac{c}{3}\right), 
\end{equation}
with the help of Fatou's lemma. The difference due to this replacement can be evaluated by using the modulus of continuity of $b$ in the variables $x, \xi$. Similarly the stochastic integrals can be tackled, which is explained in what follows (in the next steps of the proof). Denote
$$
D^n_{M,N,\epsilon} := \Omega \setminus (A^n_{M,\epsilon} \cup  B^n_{N,\epsilon}), \quad n\ge 0.
$$
Notice that $\inf_{n\ge 0}\P(D^n_{M,N,\epsilon})> 1-2\epsilon$ and that 
$$
1(D^n_{M,N,\epsilon}) \, \sup_{0\le t\le T}|\tilde Z^n_t| \vee |\tilde \zeta^n_t| \le M, \quad n\ge 0,  
$$
and 
$$
1(D^n_{M,N,\epsilon}) \, 
1\left(\bigcup_{k,j=1}^N\sup_{0\le t\le T}|\tilde X^n_t - \phi^k_t| \vee |\tilde \xi^n_t - \phi^j_t| \le \epsilon\right) = 1(D^n_{M,N,\epsilon}).
$$
Let
$$
D^{n,k,j}_{M,N,\epsilon} := \left(\omega: \sup_{0\le t\le T}|\tilde X^n_t - \phi^k_t| \vee |\tilde \xi^n_t - \phi^j_t| \le \epsilon\right).
$$

Denote for a chosen couple $(\phi^k,\phi^j)$
\begin{align*}
g^{n,n_0,k,j}(s,y,\eta):=b_0^{n}(s,\phi^k_s,y,\phi^j_s,\eta)- b_0^{n_0}(s,\phi^k_s,y,\phi^j_s,\eta), \\\\
g^{n,k,j}(s,y,\eta):=b_0^{n}(s,\phi^k_s,y,\phi^j_s,\eta)- b_0^{}(s,\phi^k_s,y,\phi^j_s,\eta).
\end{align*}
Then the first summand $I^1$  may be estimated by the BCM   inequality  as follows:
\begin{align*}
&I^1 \le \frac3{c} \, \mathsf E  \int_0^T C \,\mathsf{E}^3 |b_0^{n}(s,\tilde Z^{n}_s,\tilde  \zeta^{n}_s)- b_0^{n_0}(s,\tilde Z^n_s, \tilde \zeta^n_s)|\,ds
 \\\\
&= C \, \mathsf E  (1(D^n_{M,N,\epsilon}) + 1(A^n_{M,\epsilon} \cup  B^n_{N,\epsilon})) \int_0^T  |b_0^{n}(s,\tilde Z^{n}_s,\tilde  \zeta^{n}_s)- b_0^{n_0}(s,\tilde Z^n_s, \tilde \zeta^n_s)|\,ds.
\end{align*}
{  Due to (\ref{eps1}), (\ref{eps2}), and (\ref{linear}) we have}
\begin{align*}
\mathsf E  1(A^n_{M,\epsilon} \cup  B^n_{N,\epsilon}) \int_0^T  |b_0^{n}(s,\tilde Z^{n}_s,\tilde  \zeta^{n}_s)- b_0^{n_0}(s,\tilde Z^n_s, \tilde \zeta^n_s)|\,ds \le C\epsilon.
\end{align*}
So, it remains to evaluate the term
\begin{align*}
&\mathsf E  1(D^n_{M,N,\epsilon})\int_0^T  |b_0^{n}(s,\tilde Z^{n}_s,\tilde  \zeta^{n}_s)- b_0^{n_0}(s,\tilde Z^n_s, \tilde \zeta^n_s)|\,ds 
\\\\
&\le \sum_{k,j=1}^{N}\mathsf E  1(D^{n,k,j}_{M,N,\epsilon})\int_0^T  |b_0^{n}(s,\tilde Z^{n}_s,\tilde  \zeta^{n}_s)- b_0^{n_0}(s,\tilde Z^n_s, \tilde \zeta^n_s)|\,ds 
\end{align*}
We have for any $k,j\le N$ 
\begin{align*}
&\mathsf E  1(D^{n,k,j}_{M,N,\epsilon})\int_0^T  |b_0^{n}(s,\tilde Z^{n}_s,\tilde  \zeta^{n}_s)- b_0^{n_0}(s,\tilde Z^n_s, \tilde \zeta^n_s)|\,ds 
 \\\\
&\le \mathsf E  1(D^{n,k,j}_{M,N,\epsilon})\int_0^t  \left| b_0^{n}(s,\tilde X^{n}_s, \tilde Y^{n}_s, \tilde \xi^n_s, \tilde \eta^{n}_s)  -  b_0^{n}(s,\phi^k_s, \tilde Y^{n}_s, \phi^j_s, \tilde \eta^{n}_s)\right|ds
 \\\\
&+\mathsf E  1(D^{n,k,j}_{M,N,\epsilon})\int_0^t  \left| b_0^{n_0}(s,\tilde X^{n}_s, \tilde Y^{n}_s, \tilde \xi^n_s, \tilde \eta^{n}_s)  -  b_0^{n_0}(s,\phi^k_s, \tilde Y^{n}_s, \phi^j_s, \tilde \eta^{n}_s)\right|ds
 \\\\
&+\mathsf E  1(D^{n,k,j}_{M,N,\epsilon})\int_0^t  \left| b_0^{n}(s,\phi^k_s, \tilde Y^{n}_s, \phi^j_s,\tilde \eta^{n}_s)  -   b_0^{n_0}(s,\phi^k_s, \tilde Y^{n}_s, \phi^j_s,\tilde \eta^{n}_s)\right|ds
 \\\\
&\le C\epsilon + \mathsf E  1(D^{n,k,j}_{M,N,\epsilon})\int_0^t  \left|b_0^{n}(s,\phi^k_s, \tilde Y^{n}_s, \phi^j_s, \tilde \xi^{n}_s)  -  b_0^{n_0}(s,\phi^k_s, \tilde Y^{n}_s, \phi^j_s, \tilde \eta^{n}_s)\right|ds.
\end{align*}
By virtue of   Krylov's estimate (see  the Theorems 2.4.1 or 2.3.4 in \cite{Kry}) 
\begin{align}\label{gnn0}
&\mathsf E  1(D^{n,k,j}_{M,N,\epsilon})\int_0^t  \left|b_0^{n}(s,\phi^k_s, \tilde Y^{n}_s, \phi^j_s, \tilde \eta^{n}_s)  -  b_0^{n_0}(s,\phi^k_s, \tilde Y^{n}_s, \phi^j_s, \tilde \eta^{n}_s)\right|ds
  \nonumber\\\nonumber\\ 
&= \mathsf E  1(D^{n,k,j}_{M,N,\epsilon})\int_0^t  |g|^{n,n_0,k,j}(s,\tilde Y^{n}_s, \tilde \eta^{n}_s)ds
 \nonumber\\\\ \nonumber
&\le N_R \left(\|g^{n,k,j}\|_{L_{2d+1}([0,T]\times B_R\times B_R)} 
+ \|g^{n_0,k,j}\|_{L_{2d+1}([0,T]\times B_R\times B_R)}\right)
\to 0, 
\end{align}
as $n, n_0 \to \infty$ for each $R$, because of  the well-known property of mollified functions. Hence, overall, we obtain that
$$
I^1 \to 0, \quad n, n_0\to\infty.
$$

\medskip

Further, the second term admits the bound
\begin{align*}
&I^2 =  \mathsf P \left(\left|\int_0^t  \left( ( \mathsf{E}^3 b_0^{n_0}(s,\tilde Z^{n}_s, \tilde \zeta^{n}_s))  -  ( \mathsf{E}^3b_0^{n_0}(s,\tilde Z^0_s,\tilde  \zeta^0_s))\right)ds\right| > \frac{c}3\right)
 \\\\
&\le C \mathsf E \int_0^t  \mathsf{E}^3 \left|b_0^{n_0}(s,\tilde Z^{n}_s, \tilde \zeta^{n}_s))  -  (b_0^{n_0}(s,\tilde Z^0_s,\tilde  \zeta^0_s))\right|ds 
 \\\\
&\le C \mathsf E \int_0^t  \left|b_0^{n_0}(s,\tilde Z^{n}_s, \tilde \zeta^{n}_s))  -  (b_0^{n_0}(s,\tilde Z^0_s,\tilde  \zeta^0_s))\right|ds  \to 0, \quad n\to\infty,
\end{align*}
due to the Lebesgue bounded convergence theorem. Hence, for each $n_0$
$$
\lim_{n\to\infty} I^2 = 0, 
$$
and therefore
$$
\lim_{n_0\to\infty}\lim_{n\to\infty} I^2 = 0.
$$

\medskip

The  term $I^3$ can be considered similarly to $I^1$, { using the events $A^0_{M,\epsilon}$ and $B^0_{M,\epsilon}$ instead of $A^n_{M,\epsilon}$ and $B^n_{M,\epsilon}$}. There is just one nuance that it is not known in advance whether or not the limiting processes $\tilde Z^0, \tilde\zeta^0$ are diffusions. However, it is explained in \cite[section II.6]{Kry}; see also some details in \cite[proof of inequality (2.17)]{VerMish20}. The main point is the extension to the limiting process $(\tilde Z^0,\tilde\zeta^0)$ of Krylov's bound for diffusions $(\tilde Z^n,\tilde\zeta^n), n\ge 1$: ($1^\circ$) as a first step these bounds are proved for the limiting process in the argument of  continuous functions from $L_{2d+1}$,  and ($2^\circ$) as a second step this extension is generalised to any nonnegative Borel measurable functions using the property of regularity of probability measures in finite-dimensional Euclidean spaces. The details may be read in the cited sources. 
Hence, by the properties of the mollified functions  it follows that
$$
\lim_{n_0\to\infty}I^3 = 0.
$$
The convergence  (\ref{limNb0}) is, thus, proved.

\medskip

\noindent
{\bf 6}. Let us show for stochastic integrals in (\ref{limNs0}) that for any $c, \epsilon>0$ there exists $C>0$ such that  
\begin{equation}\label{sest0}
\mathsf P\left(\left\|\int_0^t  (\mathsf{E}^3 \sigma^n(s,\tilde Z^n_s,  \tilde \zeta^n_s)) d\tilde W^n_s  - \int_0^t   (\mathsf{E}^3 \sigma(s,\tilde Z^0_s,  \tilde \zeta^0_s)) d\tilde W^0_s\right\| > c\right)
 < C \epsilon,
\end{equation}
if $n$ is large enough. The task is similar to the convergence of Lebesgue integrals related to the coefficient $b$ studied in the previous steps of the proof. The additional obstacle is that we have to show convergence of the difference of stochastic integrals driven by different Wiener processes \(\tilde W^n \) and \(\tilde W^0\) in (\ref{sest0}) 
with {\em continuous} and bounded integrands 
$f^n_s:=  \mathsf{E}^3 \sigma^{n_0}(s,{ \tilde  Z^n_s, \tilde  \zeta^n_s})$ and $f^0_s:=  \mathsf{E}^3 \sigma^{n_0}(s,{ \tilde  Z_s^0, \tilde  \zeta_s^0})$ in Skorokhod's lemma, which does require such a boundedness {  and, at least, a continuity in probability}. Let us show the details. We have,  

\begin{align*}
&\mathsf P\left(\left\|\int_0^t  (\mathsf{E}^3 \sigma^n(s,{ \tilde Z^n_s,  \tilde \zeta^n_s})) d\tilde W^n_s  - \int_0^t   (\mathsf{E}^3 \sigma(s,{ \tilde Z_s^0,  \tilde \zeta_s^{ 0}})) d\tilde W_s\right\| > c\right)
 \\\\
&\le \mathsf P\left(\left\|\int_0^t  (\mathsf{E}^3 \sigma^n(s,{ \tilde Z^n_s,  \tilde \zeta^n_s})) d\tilde W^n_s  - \int_0^t   (\mathsf{E}^3 \sigma ^{n_0}(s,{ \tilde Z^n_s,  \tilde \zeta^n_s})) d\tilde W^n_s\right\| > c/3\right)
 \\\\
&+ \mathsf P\left(\left\|\int_0^t  (\mathsf{E}^3 \sigma^{n_0}(s,{ \tilde Z^n_s,  \tilde \zeta^n_s})) d\tilde W^n_s  - \int_0^t   (\mathsf{E}^3 \sigma ^{n_0}(s,{ \tilde Z_s^0,  \tilde \zeta_s^0})) d\tilde W_s^0\right\| > c/3\right)
 \\\\
&+\mathsf P\left(\left\|\int_0^t  (\mathsf{E}^3 \sigma^{n_0}(s,{ \tilde Z_s^0,  \tilde \zeta_s^0})) d\tilde W_s^0  - \int_0^t   (\mathsf{E}^3 \sigma ^{}(s,{ \tilde Z_s^0,  \tilde \zeta_s^0})) d\tilde W_s^0\right\| > c/3\right)
 \\\\
&\le C \mathsf E \left\|\int_0^t  \mathsf{E}^3 (\sigma^n(s,{ \tilde Z^n_s,  \tilde \zeta^n_s}))  -    \sigma ^{n_0}(s,{ \tilde Z^n_s,  \tilde \zeta^n_s})) d\tilde W^n_s\right\|^2
 \\\\
&+C \mathsf E \left\|\int_0^t  (\mathsf{E}^3 \sigma^{n_0}(s,{ \tilde Z^n_s,  \tilde \zeta^n_s})) d\tilde W^n_s  - \int_0^t   (\mathsf{E}^3 \sigma ^{n_0}(s,{ \tilde Z_s^0,  \tilde \zeta_s^0})) d\tilde W_s^0\right\| 
 \\\\
&+ C \mathsf E \left\|\int_0^t  \mathsf{E}^3 (\sigma^{n_0}(s,{ \tilde Z_s^0,  \tilde \zeta_s^0}))   -     \sigma ^{}(s,{ \tilde Z_s^0,  \tilde \zeta_s^0})) d\tilde W_s^0\right\|^2
 \\\\
&\le C \mathsf E \int_0^t  \mathsf{E}^3 \left\|(\sigma^n(s,{ \tilde Z^n_s,  \tilde \zeta^n_s}))  -    \sigma ^{n_0}(s,{ \tilde Z^n_s,  \tilde \zeta^n_s}))\right\|^2 ds
 \\\\
&+C \mathsf E \left\|\int_0^t  \underbrace{\mathsf{E}^3 \sigma^{n_0}(s,\tilde Z^n_s,  \tilde \zeta^n_s)}_{:=f^n(s,\omega)} d\tilde W^n_s  - \int_0^t   \underbrace{\mathsf{E}^3 \sigma ^{n_0}(s,\tilde Z^0_s,  \tilde \zeta^0_s)}_{:=f^0(s,\omega)} d\tilde W^0_s\right\| 
 \\\\
&+ C \mathsf E \int_0^t  \mathsf{E}^3 \left\|(\sigma^{n_0}(s,{ \tilde Z_s^0,  \tilde \zeta_s^0}))   -     \sigma ^{}(s,{ \tilde Z_s^0,  \tilde \zeta_s^0}))\right\|^2 ds
 \\\\
&=: J^1 + J^2 + J^3.
\end{align*}
The terms $J^1$ and $J^3$ are tackled similarly to $I^1$ and $I^3$ from the previous steps. The additional difficulty discussed above relates to the term $J^2$. 
Skorokhod's lemma (see \cite[Lemma 4 in the Appendix]{VerMish20}) is applicable if
for {\em bounded and (stochastically) continuous  in $s$} integrands $f^n:=\mathsf{E}^3 \sigma^{n_0}(s,\tilde Z^n_s,  \tilde \zeta^n_s)$  uniformly with respect to $n\ge 0$ it holds that 
\begin{equation}\label{fn0}
f^n(s,\omega) \stackrel{\P}\to f^0(s,\omega), \quad \text{a.e. $s\le T$}\quad n\to\infty.
\end{equation}
The fact that $f^n$ is uniformly bounded follows straightforwardly from the boundedness of the function $\sigma$. Further, $f^n(s,\omega)$ is continuous in $s$ a.s. uniformly in $n\ge 0$ because of the continuity of $\sigma^{n_0}$ in all variables and due to the uniform stochastic continuity of all processes ${ \tilde Z^n_s}$ for $n\ge 0$ (see (\ref{kol2t})). Finally, the convergence (\ref{fn0}) in probability {\em for all $s\le T$} follows from the following little calculus:
\begin{align*}
&\mathsf E \|f^n(s,\omega) - f^0(s,\omega)\| 
= \mathsf E \|\mathsf{E}^3 \sigma^{n_0}(s,\tilde Z^n_s,  \tilde \zeta^n_s) - \mathsf{E}^3 \sigma^{n_0}(s,\tilde Z^0_s,  \tilde \zeta^0_s)\| 
 \\\\
&= \mathsf E \|\mathsf{E}^3 \left(\sigma^{n_0}(s,\tilde Z^n_s,  \tilde \zeta^n_s) -  \sigma^{n_0}(s,\tilde Z^0_s,  \tilde \zeta^0_s) \right)\| 
\le \mathsf E \mathsf{E}^3 \|\sigma^{n_0}(s,\tilde Z^n_s,  \tilde \zeta^n_s) - \mathsf{E}^3 \sigma^{n_0}(s,\tilde Z^0_s,  \tilde \zeta^0_s)\| 
 \\\\
&= \mathsf E \|\sigma^{n_0}(s,\tilde Z^n_s,  \tilde \zeta^n_s) - \mathsf{E}^3 \sigma^{n_0}(s,\tilde Z^0_s,  \tilde \zeta^0_s)\|  \to 0, \quad n\to\infty.
\end{align*}
The latter convergence holds true by virtue of Lebesgue's bounded convergence theorem.
Hence, in this way we obtain that
$$
 J^1 + J^2 + J^3 \to 0, \quad n, n_0\to\infty.
$$
So, the convergence (\ref{limNs0}) holds true, which  along with (\ref{limNb0})  leads to the equation (\ref{limitsde}) and, hence, completes the proof of the theorem. \hfill QED


\end{document}